\documentclass{amsart}
\usepackage{amscd,amssymb,amsmath,verbatim,hyperref,
epsfig,latexsym,subfigure}
\usepackage[arrow,matrix,graph,frame,poly,arc,tips]{xy}
 
\begin{document}

\newcommand{\mmbox}[1]{\mbox{${#1}$}}
\newcommand{\proj}[1]{\mmbox{{\bf P}^{#1}}}
\newcommand{\affine}[1]{\mmbox{{\bf A}^{#1}}}
\newcommand{\Ann}[1]{\mmbox{{\rm Ann}({#1})}}
\newcommand{\caps}[3]{\mmbox{{#1}_{#2} \cap \ldots \cap {#1}_{#3}}}

\newcommand{\Tor}{\mathop{\rm Tor}\nolimits}
\newcommand{\Ext}{\mathop{\rm Ext}\nolimits}
\newcommand{\Hom}{\mathop{\rm Hom}\nolimits}
\newcommand{\im}{\mathop{\rm Im}\nolimits}
\newcommand{\rank}{\mathop{\rm rank}\nolimits}
\newcommand{\supp}{\mathop{\rm supp}\nolimits}

\newcommand{\arrow}[1]{\stackrel{#1}{\longrightarrow}}
\sloppy

\newtheorem{defn0}{Definition}[section] 
\newtheorem{prop0}[defn0]{Proposition}
\newtheorem{conj0}[defn0]{Conjecture}
\newtheorem{thm0}[defn0]{Theorem}
\newtheorem{lemma0}[defn0]{Lemma}
\newtheorem{corollary0}[defn0]{Corollary}
\newtheorem{example0}[defn0]{Example}

\newenvironment{defn}{\begin{defn0}}{\end{defn0}}
\newenvironment{prop}{\begin{prop0}}{\end{prop0}}
\newenvironment{conj}{\begin{conj0}}{\end{conj0}}
\newenvironment{thm}{\begin{thm0}}{\end{thm0}}
\newenvironment{lemma}{\begin{lemma0}}{\end{lemma0}}
\newenvironment{corollary}{\begin{corollary0}}{\end{corollary0}}
\newenvironment{example}{\begin{example0}\rm}{\end{example0}}

\newcommand{\defref}[1]{Definition~\ref{#1}}
\newcommand{\propref}[1]{Proposition~\ref{#1}}
\newcommand{\thmref}[1]{Theorem~\ref{#1}}
\newcommand{\lemref}[1]{Lemma~\ref{#1}}
\newcommand{\corref}[1]{Corollary~\ref{#1}}
\newcommand{\exref}[1]{Example~\ref{#1}}
\newcommand{\secref}[1]{Section~\ref{#1}}

\title[Linear systems on a special rational surface]%
{Linear systems on a special rational surface}
 
\author[Henry K. Schenck]{Henry K. Schenck}
\address{Department of Mathematics,
Texas A\&M University, College Station, TX 77843}
\email{\href{mailto:schenck@math.tamu.edu}{schenck@math.tamu.edu}}
\urladdr{\href{http://www.math.tamu.edu/~schenck/}%
{http://www.math.tamu.edu/\~{}schenck}}
 
\thanks{partially supported by NSF grant DMS03-11142 and NSA grant MDA904-0301-0006.}
 
\begin{abstract}
We study the Hilbert series of two families of ideals generated 
by powers of linear forms in $\mathbb{K}[x_1,x_2,x_3]$. Using the 
results of Emsalem-Iarrobino, we formulate this as a question about 
fatpoints in $\mathbb{P}^{2}$. This is equivalent to studying the 
dimension of a linear system on a blow up of $\mathbb{P}^2$. 
We determine the classes of the negative curves, then apply an
algorithm of Harbourne to reduce to an effective, nef divisor.
Combining Harbourne's results on rational surfaces with $K^2 > 0$ and
Riemann-Roch yields a formula for the Hilbert series.
For one family of ideals, this proves the $n=3$ case of a conjecture 
posed by Postnikov and Shapiro ``as a challenge to the 
commutative algebra community'' (after this proof was communicated
to them, they found a combinatorial proof for all
$n$). For the second family of ideals, it yields a formula which 
Postnikov and Shapiro were unable to obtain via combinatorial 
methods. We conjecture a formula for the minimal free resolution 
of one family of ideals, and show that a member of the second family 
of ideals provides a counterexample to a conjecture made by 
Postnikov and Shapiro in \cite{ps}.
  \end{abstract}    
\maketitle
\section{Introduction}
\setcounter{defn0}{0}
Let $R=\mathbb{K}[x_1,\ldots,x_n]$ be a polynomial ring. 
We consider the following families of $2^n-1$ generated ideals:
$$I_\phi=\langle x_1^{\phi(1)},\ldots, x_n^{\phi(1)}, (x_1x_2)^{\phi(2)},\ldots,
(x_{i_1}\cdots x_{i_r})^{\phi(r)}\ldots \rangle$$
$$J_\phi=\langle x_1^{\phi(1)},\ldots, x_n^{\phi(1)}, (x_1+x_2)^{2 \phi(2)},\ldots,
(x_{i_1}+\cdots +x_{i_r})^{r \phi(r)}\ldots \rangle,$$
where $\phi$ is either a {\em linear degree function}:
 $$\phi(r)=l+k(n-r)>0,\mbox{ } k,l \in \mathbb{N}\mbox{ (the positive integers)}.$$
or an {\em almost linear degree function}: $\exists k \in \mathbb{N}$ 
such that 
$$\phi(r)-\phi(r+1)=k\mbox{ or }k+1\mbox{ }\forall r \in \{1,\ldots,n-1\}.$$

In \cite{ps1}, Postnikov and Shapiro conjecture that for a linear 
degree function $\phi$, the Hilbert series $P(R/I_\phi,t)$ is equal 
to the Hilbert series $P(R/J_\phi, t).$ They gave
a proof for $n=2$, and for any $n$, when $\phi(r) = n+1-r$.
If $\phi$ is almost linear, they observe that the Hilbert series 
are also often equal, although not always, as there are
counterexamples when $n=4$. In this paper, we prove:
\begin{thm}
If $\phi$ is linear or almost linear, $n=3$, and 
$char(\mathbb{K})=0$ or $char(\mathbb{K}) > \sum_{i=1}^3\phi(i)-2$,
then $P(R/I_\phi,t) = P(R/J_\phi, t).$
\end{thm}
Theorem 1.1 is proved via algebraic geometry. While none of the 
previous work on $P(R/J_\phi, t)$ (\cite{ps1}, \cite{ps}, \cite{ssp})
made any assumptions on the characteristic of $\mathbb{K}$, some conditions 
are in fact necessary. We illustrate this in Example 3.2, which 
shows that in low positive characteristic the conclusion of 
Theorem 1.1 can fail in very simple cases, and for very simple reasons.

After the proof of Theorem 1.1 was communicated to them \cite{ssp}, 
Postnikov and Shapiro \cite{ps} were able to use combinatorial 
techniques (the correspondence between parking functions and labeled
trees) to give a proof that if $\phi$ is linear then  
$P(R/I_\phi,t) = P(R/J_\phi, t)$ for any $n$. However, their
methods do not handle the almost linear case.
They note that except in very special cases a 
Gr\"obner basis approach {\it cannot} work because the
monomial generators do not lie on the boundary of the 
Newton polytopes of their polynomial deformations.
Their motivation for the work was an effort to
generalize earlier results \cite{pss} on the algebra generated by
curvature forms on the generalized flag manifold.
 
For a monomial ideal $I$, the Taylor resolution (\cite{t}, or 
\cite{eis}, Exercise 17.11) is a generalization of the Koszul 
complex which gives an
explicit free (generally non-minimal) resolution of $I$.
Thus, $P(R/I_\phi, t)$ is known; the interesting case is 
$R/J_\phi$. For a linear degree function, it seems possible
that $J_\phi$ and $I_\phi$ have equal graded Betti numbers; this
is consistent with Theorem 1.1. We give an example of
an almost linear degree function where $P(R/I_\phi, t)=P(R/J_\phi,
t)$, but where the graded Betti numbers differ; this shows that
Conjecture 6.10 in \cite{ps} is false.

We first use Macaulay's inverse 
systems approach to relate the dimension of $J_\phi$ in
degree $j$ to the degree $j$ piece of an ideal $F$ of fatpoints 
on $\mathbb{P}^2$. $J_\phi$ is generated by seven forms, so 
$F$ is supported at seven points. Following Harbourne,
we blow up $\mathbb{P}^2$ at the seven points and 
study a linear system on the resulting surface $X$; because only seven
points were blown up we can determine all the negative classes on $X$,
which allows us to determine $P(R/J_\phi, t)$. When $n=2$, $J_\phi$ 
behaves as if
the forms were generic. In \cite{f}, Fr\"oberg made a conjecture
about the behavior of the Hilbert series of an ideal generated
by generic forms in $n$ variables, and proved the conjecture when $n=2$;
the $n=3$ case was subsequently solved by Anick in \cite{a}.
We show that the deviation of the dimension of $(J_\phi)_j$ from the dimension
of an ideal generated by generic forms of the same degrees
is measured by $H^1(D_j)$, where $D_j$ is a divisor on $X$ corresponding
to the degree $j$ piece of a fatpoint ideal on  $\mathbb{P}^2$.\newline
\newline
\noindent {\bf Acknowledgments} I thank Boris Shapiro for explaining
the problem to me, M.F.I.-Oberwolfach for bringing us together, 
and Brian Harbourne and Tony Iarrobino 
for several enlightening conversations. The Macaulay 2 software
package was also a help, as were comments from an anonymous
referee. 

\section{Linear Forms and Fatpoints}
In \cite{ei}, Emsalem and Iarrobino proved that there is 
a close connection between ideals generated by powers of linear
forms, and ideals of fatpoints. 
Let $p_i = [p_{i0}:p_{i1}:\cdots :p_{in}] \in \mathbb{P}^n$, $I(p_i) = 
\wp_i \subseteq R = \mathbb{K}[x_0,\ldots, x_n]$, and $L_{p_i} = 
\sum_{j=0}^np_{i_j}y_j$. Let $\{ p_1, \ldots , p_m \} \subseteq
\mathbb{P}^n$ be a set of distinct points. 
A {\em fatpoints ideal} is an ideal of the 
form $$F = \cap_{i=1}^m \wp_i^{\alpha_i}.$$ Let
$S = \mathbb{K}[y_0,\ldots, y_n]$, and define an action of $R$ on $S$ by 
partial differentiation, i.e. $x_j \cdot y_i = \partial(y_i)/\partial(y_j)$.
We think of $S$ both as a ring, isomorphic
to $R$, and as an $R$-module. Since $F$ is a submodule of
$R$, it acts on $S$, and we can ask what elements of $S$ are annihilated
by this action. The set of such elements is denoted by $F^{-1}$. 
The essential result of Emsalem and Iarrobino is that for $j \gg 0$, $(F^{-1})_j = 
\langle L_{p_1}^{j-\alpha_1+1}, \ldots, L_{p_m}^{j-\alpha_m+1} \rangle_j$, and that
$\dim_\mathbb{K}(F^{-1})_j = \dim_\mathbb{K}(R/F)_j$. This generalizes the classical Terracini
lemma \cite{t1}, \cite{t2}, where the $\alpha_i$ are all two.

\begin{thm} $($Emsalem and Iarrobino, \cite{ei}$)$ Let
$F$ be an ideal of fatpoints: $$F = \wp_1^{n_1+1} \cap \cdots
\cap \wp_m^{n_m+1}.$$ If $char(\mathbb{K})=0$ or $char(\mathbb{K})>j$, then 

$$(F^{-1})_j = \begin{cases} S_j & \mbox{for } j \leq \max\{n_i\} \cr
                               &                         \cr
       L_{p_1}^{j-n_1}S_{n_1} + \cdots + L_{p_m}^{j-n_m}S_{n_m} &
                               \mbox{for }
j\geq\max\{ n_i + 1\} \end{cases}           
$$
and  
$$
\dim _\mathbb{K}(F^{-1})_j = \dim _\mathbb{K}(R/F)_j.
$$
\end{thm}
Suppose we have an ideal generated by powers of linear forms, and 
for each $j \in \mathbb{N}$, we wish to compute the dimension of 
$$
 \langle L_{p_1}^{t_1}, \cdots , L_{p_m}^{t_m} \rangle_j. 
$$
Since the $t_i$ are fixed, to apply the approach above we fix a degree $j$. 
Put
$$F(j) = \wp_1^{j-t_1+1} \cap \cdots \cap \wp_m^{j-t_m+1}.$$ 
Then 
$$
\mbox{dim}_\mathbb{K} \langle L_{p_1}^{t_1}, \cdots , L_{p_m}^{t_m} \rangle_j
= \mbox{dim}_\mathbb{K}(R/F(j))_j.$$
Hence, we will be studying an infinite family of ideals of fatpoints. 
For additional information on this correspondence we refer to 
the original paper \cite{ei}, Geramita \cite{g}, or Macaulay \cite{m};
an improved version of Theorem 2.1 is given in \cite{ik}. 
In the next section we apply Theorem 2.1 to study fatpoints ideals 
in $\mathbb{P}^2$ related to linear degree functions, so henceforth 
we require that the characteristic of our field be zero or greater
than $\phi(1)+\phi(2)+\phi(3)-2$.

\section{Blow-ups, Rational Surfaces, and the problem}
Recall now the problem: given a linear degree function $$\phi(r)=l+k(n-r)>0,\mbox{ } k,l \in \mathbb{N},$$
determine the Hilbert series of the quotient of $R$ by the
ideal generated by (over all tuples) 
$$(x_{i_1}+\cdots+x_{i_r})^{r \phi(r)}.$$
As we saw in the previous section, this is equivalent to
determining the Hilbert function of a {\it family} of ideals of
fatpoints. In \cite{n} Nagata studied fatpoints in $\mathbb{P}^2$,
obtaining results for small numbers of points and posing a number 
of conjectures. This continues to be a very active area of research;
see \cite{h1} or \cite{m1} for overviews. Surprisingly, even 
for a set of points in general position in $\mathbb{P}^2$ the 
Hilbert function is unknown. A formula is conjectured by
Segre, Harbourne, Gimigliano, and Hirschowitz, and there has been
substantial recent progress: for more, see e.g. \cite{bz},
\cite{cm}, \cite{hr}. For the remainder of the paper, we 
specialize to the case of 
three variables. Given $\phi$, we want to find the Hilbert 
function of the ideal
$$\langle x^{\phi(1)}, y^{\phi(1)},z^{\phi(1)}, (x+y)^{2\phi(2)}, (x+z)^{2\phi(2)},(y+z)^{2\phi(2)}, (x+y+z)^{3\phi(3)} \rangle.$$
Consider the following seven points of $\mathbb{P}^2$:
\begin{center}
\epsfig{file=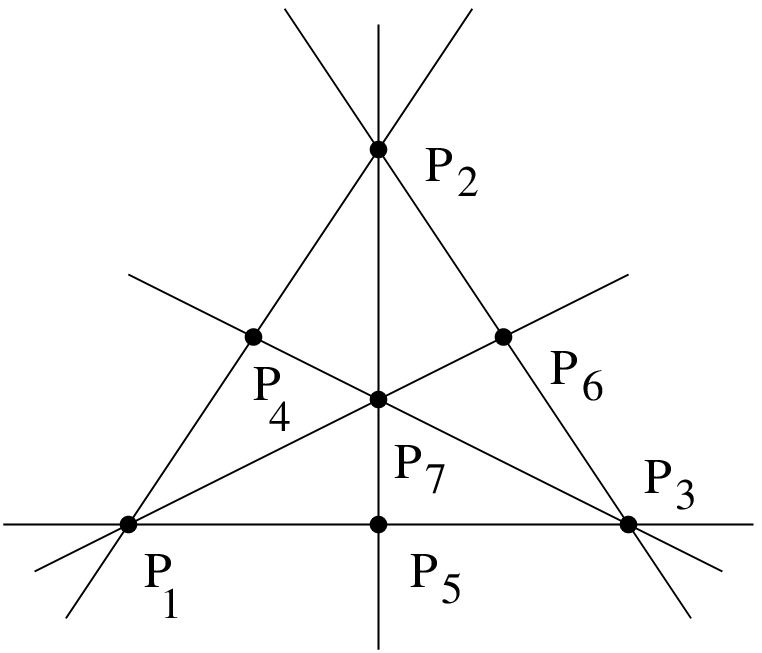,height=1.7in,width=1.5in}
\end{center}
\begin{center}
$\begin{array}{ccc}
p_1&= &(1:0:0) \\
p_2&= &(0:1:0) \\
p_3&= &(0:0:1) \\
p_4&= &(1:1:0) \\
p_5&= &(1:0:1) \\
p_6&= &(0:1:1) \\
p_7&= &(1:1:1)
\end{array}$
\end{center}
and let $\wp_i=I(p_i)$. Define $(a)_+ = max\{a,0\}$. 
The results of the previous section show
that $\mbox{dim}_\mathbb{K}(J_\phi)_j = \mbox{dim}_\mathbb{K}(R/F(j))_j$, where 
$$F(j)=\langle \wp_1^{a_1} \cap \wp_2^{a_2} \cap \wp_3^{a_3} \cap \wp_4^{a_4}\cap \wp_5^{a_5}\cap \wp_6^{a_6} \cap \wp_7^{a_7}\rangle,$$
and $a_1=a_2=a_3=(j-\phi(1)+1)_+$, $a_4=a_5=a_6=(j-2\phi(2)+1)_+$, and $a_7=(j-3\phi(3)+1)_+$. The key to solving this problem is work of Harbourne \cite{h2} 
which shows how to determine the dimension of a 
linear system on a blow up of $\mathbb{P}^2$ at eight or fewer points.

We begin by recalling that there is a correspondence between the 
graded pieces of an ideal of fatpoints $F$ and the global sections of a
certain line bundle on the surface $X$ which is the blow up of $\mathbb{P}^2$
at the points. Let $E_i$ be the class of the exceptional curve
over the point $p_i$, and $E_0$ the pullback of a line on $\mathbb{P}^2$.
Put $$D_j = jE_0-\sum\limits_{i=1}^7 a_iE_i,$$
with $a_i$ as above. The canonical divisor on $X$ is:
$$K_X = -3E_0+\sum\limits_{i=1}^7 E_i.$$
Since $j>0$, $h^2(D_j)=0$ and Riemann-Roch yields:
$$\chi(D_j)= h^0(D_j)-h^1(D_j) =  \frac{D_j^2-D_j K_X}{2}+1.$$ 
We compute
$$\frac{D_j^2-D_j K_X}{2}+1 = {j+2 \choose 2} -3{j+2-\phi(1) \choose 2}-
3{j+2-2\phi(2) \choose 2}-{j+2-3\phi(3) \choose 2}.$$
Since $$F(j)_j=H^0(D_j),$$
the dimension of $(J_\phi)_j$ is given by:
$$3{j+2-\phi(1) \choose 2} + 3{j+2-2\phi(2) \choose 2}+{j+2-3\phi(3) \choose 2}-h^1(D_j).$$
Notice the connection to Anick's result. $J_\phi$ is generated by three forms of degree $\phi(1)$, three forms of degree $2\phi(2)$,
and one form of degree $3\phi(3)$. Thus, the maximal dimension of $\langle J_\phi \rangle_j$ is obtained by multiplying all monomials of degree $j-\phi(1)$ with
the generators of degree $\phi(1)$, and similarly for the other generators, so the maximal dimension of $\langle J_\phi \rangle_j$ is
$$3{j-\phi(1)+2 \choose 2}+ 3{j-2\phi(2)+2 \choose 2}+{j-3\phi(3)+2 \choose 2}.$$
In other words, $(J_\phi)_j$ is as large as possible iff $h^1(D_j) =  0$, i.e.
exactly when the divisor $D_j$ is nonspecial. The original problem may
be restated as follows: for a linear degree function, determine the number of global sections i.e. compute $h^0(D_j).$  
\begin{example} For $n=3$, take $l=3$ and $k=1$.
Let $G_\phi$ be generated by generic forms of the same degree as the
generators of $J_\phi$. If $char(\mathbb{K})=0$ or
$char(\mathbb{K}) \ge 7$, then
$$P(R/G_\phi,t)=1+ 3t+ 6t^2+ 10t^3+ 15t^4+ 18t^5+ 19t^6+ 18t^7+ 12t^8$$
$$P(R/I_\phi,t)=P(R/J_\phi,t)=1+ 3t+ 6t^2+ 10t^3+ 15t^4+ 18t^5+ 19t^6+ 18t^7+ 12t^8+6t^9.$$
There are many examples of this type, which is to be
expected, since monomial ideals will very rarely
have generic Hilbert series.
\end{example}

\begin{example} For $n=3$, take $l=1$ and $k=1$. We have that
  $\phi(1)=3, \phi(2)=2, \phi(3)=1$, so 
$$I_\phi = \langle x^3, y^3,z^3, (xy)^2, (yz)^2,(xz)^2, xyz \rangle.$$
$$J_\phi = \langle x^3, y^3,z^3, (x+y)^4, (y+z)^4,(x+z)^4, (x+y+z)^3 \rangle.$$
If $char(\mathbb{K}) =2$, 
then it is obvious that $J_\phi$ is in fact minimally
generated by $\langle x^3, y^3,z^3, (x+y+z)^3 \rangle.$ A computation 
shows that in this case
$$P(R/I_\phi,t)=1+ 3t+ 6t^2+ 6t^3, \mbox{ and }P(R/J_\phi,t)=1+ 3t+ 6t^2+ 6t^3+ 4t^4+ t^5.$$
\end{example}

\begin{example} For $n=3$, take $l=1$ and $k=2$. If
  $char(\mathbb{K})=0$ then $J_\phi$ is generated by
$$\langle x^5, y^5,z^5, (x+y)^6, (y+z)^6,(x+z)^6, (x+y+z)^3 \rangle.$$
\begin{center}
$\begin{array}{ccccc}
j & & \chi(D_j)&  &\mbox{dim}_\mathbb{K}(R/J_\phi)_j  \\
3 & & 9 & & 9 \\
4 & & 12 & & 12 \\
5 & & 12 & & 12 \\
6 & & 6 & &6\\
7 & & -6& &0
\end{array}$
\end{center}
In degree $3$, where $J_\phi$ first has a generator, the corresponding fatpoint
ideal is the ideal of the point $(1:1:1) = \langle x-z,y-z \rangle 
= F(3)$ and the dimension of $R/F(3)$ in degree $3$ is obviously $1$. 
In degree six we consider the fatpoint ideal
$$F(6)=\langle \wp_1^{2} \cap \wp_2^{2} \cap \wp_3^{2} \cap \wp_4 \cap \wp_5 \cap \wp_6 \cap \wp_7^{4}\rangle.$$
To determine the dimension of $(R/F(6))_6$, we need to understand effective
divisors on $X$. We tackle this in the next section. 
\end{example}
\section{Effective Divisors and the linear case}
We first recall some basic terminology, referring to \cite{hart} for 
additional details. Let $X$ be a smooth surface. A
divisor $D$ on $X$ is a finite integral combination of 
irreducible curves on $X$. $D$ is said to be {\it effective} if
$h^0(D)>0$, and {\it numerically effective $($nef$)$} if $D \cdot C \ge 0$ for
every effective divisor $C$, where $\cdot$ denotes the intersection
pairing. In the last section, we saw that Riemann-Roch yields a 
simple numerical formula for $h^0(D_j)-h^1(D_j)$; our goal
is to determine $h^0(D_j)$. A {\it Zariski decomposition} of an 
effective divisor $F$ is a representation $F \sim G+Z$, with 
$G$ effective and nef, $Z$ effective, and such that $h^0(F) = h^0(G).$ 
If $X$ is a blow-up of $\mathbb{P}^2$ at a small number of points, 
then Harbourne gives an algorithm to find a Zariski decomposition 
in \S 2 of \cite{h}. 

For a smooth projective rational surface $Y$ with $K_Y^2 >0$, 
Harbourne shows in Theorem 8 of \cite{h2} that if $G$ is nef, 
then $G$ is effective and $h^1(G)=h^2(G)=0$. Now let $X$ be
a blow up of $\mathbb{P}^2$ at seven points. Since $K_X^2=2$, 
finding a Zariski decomposition $D_j \sim G+Z$ and applying 
Riemann-Roch to $G$ will yield $h^0(D_j)$. To obtain a Zariski
decomposition of $D_j$ in the special 
case of a blow up at $n\le7$ points, test the divisor class
$D=D_j$ against classes $C$ of {\it negative curves}: reduced and 
irreducible curves with negative self-intersection (in this 
case, these curves generate the monoid EFF(X) of effective divisor 
classes, see \cite{h3}, Remark III.13). Subtract off any such 
class $C$ meeting $D_j$ negatively and continue
with $D=D_j-C$ either until no negative curves meet $D$ negatively,
or until $D\cdot E_0<0$. One of these two outcomes is guaranteed
to occur. If the former happens, $D_j$ was effective
and we take $G=D$; if the latter, then $D_j$ was not effective
to begin with. So the first task is to determine the negative curves.

For seven {\it general} points in $\mathbb{P}^2$, 
any negative curve is a $-1$ curve. These consist of the blow ups of the points, 
the proper transforms of lines through pairs of points, 
proper transforms of conics through any five points, and proper
transforms of cubics through six points, 
with a double point at the seventh (see \cite{man}, Theorem 26.2). 
As the points specialize, these $-1$ curves can
become reducible, but the classes of their irreducible components 
generate EFF(X). As in the previous 
section, let $E_i$ be the exceptional divisor over a point $p_i$, 
and $E_0$ the class of the pullback of a line.

\begin{lemma}
The negative curves on $X$ are the $E_i$, $i\ne 0$,
the six $-2$ curves corresponding to lines through three 
collinear points, and the three $-1$ curves corresponding
to lines through the pairs $\{p_4,p_5\},\{p_4,p_6\},\{p_5,p_6\}$:
\begin{center}
$\begin{array}{ccc}
C_{124} & = & E_0-E_1-E_2-E_4 \\
C_{135} & = &E_0-E_1-E_3-E_5 \\
C_{236} & = &E_0-E_2-E_3-E_6 \\
C_{167} & = &E_0-E_1-E_6-E_7 \\
C_{257} & = &E_0-E_2-E_5-E_7 \\
C_{347} & = &E_0-E_3-E_4-E_7 \\
C_{45} & = &E_0-E_4-E_5 \\
C_{46} & = &E_0-E_4-E_6 \\
C_{56} & = &E_0-E_5-E_6
\end{array}$
\end{center}
\end{lemma}
\begin{proof}
First, observe that $-K_X$ is 
nef. This follows since if a divisor 
$D \sim \sum_{a_i \ge 0} a_iC_i$ with $C_i$ effective and 
irreducible and $D \cdot C_i \ge 0$, then $D$ is nef: if not, then
$D\cdot F < 0 $ for some effective curve $F$, so $D\cdot F' <0$ for 
some irreducible component $F'$ of $F$. So $F'$ meets a summand of $D$
negatively; since the summand is irreducible, it must be a multiple
of $F'$. This means a summand of $D$ meets $D$ negatively, a
contradiction. In particular, $-K_X$ is nef since 
\[
-K_X = C_{167}+C_{236}+C_{45}+E_6.
\]
Any curve $C$ must satisfy the genus
formula $C^2 + C\cdot K_X = 2g_C -2.$ Since $-K_X$ is nef, 
$-K_X \cdot C \ge 0$; and $g_C$ is also 
non-negative, so the possibilities for negative curves are
$C^2=-1=K_X\cdot C$ and $C^2=-2$, $K_X\cdot C=0$. In both cases, 
there are only finitely many solutions (see \cite{man}, 25.5.3, 26.1);
we want those which are not forced by Bezout to be reducible, which
are those listed above.
\end{proof}
We return to Example $3.3$: $$D_6=6E_0-2E_1-2E_2-2E_3-E_4-E_5-E_6-4E_7.$$
We compute that $D_6 \cdot C_{167} < 0$, $(D_6-C_{167})\cdot
C_{257}<0$, 
and $(D_6-C_{167}-C_{257}) \cdot C_{347}<0$. Removing 
$C_{167}+C_{257}+ C_{347}$ from $D_6$ yields an effective, nef divisor: 
$$G=3E_0-E_1-E_2-E_3-E_7.$$ The rest is easy:
$G^2=5, GK=-5$, so $h^0(G) = 6$. Similar computations
show that this $J_\phi$ behaves generically. 

\noindent This example illustrates the crucial idea for the general case. 
Observe that the classes $C_{ijk}$ corresponding to lines with 
three collinear points are orthogonal. 
$D_j = jE_0-\sum\limits_{i=1}^7 a_iE_i$ meets
$C_{124}$ negatively iff $$j<a_1+a_2+a_4,$$ and since
$C_{124},C_{135},C_{236}$ each intersect $D_j$ with the same
multiplicity, anytime that we remove one class, 
we need to remove all three classes. Of course, this depends 
on the fact that we began with $a_1=a_2=a_3$ and 
$a_4=a_5=a_6$. In the same vein, 
$$D_j \cdot C_{167} < 0 \mbox{ iff } j<a_1+a_6+a_7,$$ 
and the curves $C_{167},C_{257},C_{347}$ 
meet $D_j$ with the same multiplicity, so again, when we remove one of 
these three, we have to remove all of them. The point is
that both of these reductions preserve the equality of
the coefficients of $\{ E_1,E_2,E_3 \}$ and of $\{ E_4,E_5,E_6 \}$. 
For the rest of the section, let $\phi$ be a linear degree
function:
\[
\begin{array}{ccc}
\phi(1) &=& l+2k\\
\phi(2) &=& l+k\\
\phi(3) &=& l
\end{array}
\] 
for some positive integers $l,k$; recall the notation $(a)_+ = max\{a,0\}$.
\begin{lemma}
$D_j$ meets $C_{124},C_{135},C_{236}$
negatively iff $j \ge \phi(1)+\phi(2)-1$. Put 
$$t_1=(j-2l-3k+2)_+ = (j-\phi(1)-\phi(2)+2)_+$$
Then $$D_j'=(j-3t_1)E_0-\sum_{i=1}^3(a_i-2t_1)E_i-\sum_{i=4}^6(a_i-t_1)E_i-a_7E_7$$
meets $C_{124},C_{135},C_{236}$ non-negatively, and has $h^0(D_j)=h^0(D_j')$.
\end{lemma}
\begin{proof}
\[
\begin{array}{ccc}
D_j \cdot C_{124} &=& j-2(j-\phi(1)+1)_+ -(j-2\phi(2)+1)_+\\
                  & = & j-2(j-l-2k+1)_+ -(j-2l-2k+1)_+
\end{array}
\]

\[
\begin{array}{ccccc}
j \le l+2k-1 & \Rightarrow  & D_j \cdot C_{124} &=& j \ge 0\\
l+2k-1 < j \le 2l+2k-1 & \Rightarrow  & D_j \cdot C_{124} &=&
-j+2l+4k-2 \ge 0\\
2l+2k-1 < j & \Rightarrow  & D_j \cdot C_{124} &=&
-2j+4l+6k-3
\end{array}
\]
So, $D_j \cdot C_{124} <0$ iff $2j \ge 4l+6k-2$ iff $j \ge 2l+3k-1 =
\phi(1)+\phi(2)-1$, and removing $t_1$ copies of 
$C_{124}+C_{135}+C_{236}$ from $D_j$ yields a divisor $D'_j$ 
which meets the classes $C_{124},C_{135},C_{236}$ non-negatively, and
has $h^0(D_j)=h^0(D_j')$.
\end{proof}
\begin{lemma}
$D_j$ meets $C_{167},C_{257},C_{347}$
negatively iff $j \ge 2\phi(2)+\phi(3)-1$. Put 
$$t_2=(j-3l-2k+2)_+=(j-2\phi(2)-\phi(3)+2)_+$$
Then $$D_j''=(j-3t_2)E_0-\sum_{i=1}^6(a_i-t_2)E_i-(a_7-3t_2)E_7$$
meets $C_{167},C_{257},C_{347}$ non-negatively, and has $h^0(D_j)=h^0(D_j'')$.
\end{lemma}
\begin{proof}
\[
\begin{array}{ccc}
D_j \cdot C_{167} &= & j-(j-\phi(1)+1)_+ -(j-2\phi(2)+1)_+ -(j-3\phi(3)+1)_+\\
 &= & j-(j-l-2k+1)_+ -(j-2l-2k+1)_+ -(j-3l+1)_+
\end{array}
\]
We have to analyze three possible cases:
\[
\begin{array}{ccccc}
l+2k-1 &\le& 2l+2k-1& \le& 3l-1\\
l+2k-1 &\le& 3l -1 &\le & 2l+2k-1\\
3l-1 &\le & l+2k-1 &\le& 2l+2k-1
\end{array}
\]
In the first case, we find 
\[
\begin{array}{ccccc}
0 \le j \le l+2k-1 & \Rightarrow  & D_j \cdot C_{167} &=& j \ge 0\\
l+2k-1< j \le 2l+2k-1 & \Rightarrow  & D_j \cdot C_{167} &=&
l+2k-1 \ge 0\\
2l+2k-1 < j \le 3l-1 & \Rightarrow  & D_j \cdot C_{167} &=&
-j+3l+4k-2 \ge 0\\
 3l-1 < j   & \Rightarrow  & D_j \cdot C_{167} &=&
-2j+6l+4k-3 
\end{array}
\]
The second case differs for the values of $j$ between $l+2k-1$ 
and $2l+2k-1$:
\[
\begin{array}{ccccc}
l+2k-1< j \le 3l-1 & \Rightarrow  & D_j \cdot C_{167} &=&
l+2k-1 \ge 0\\
3l-1< j \le 2l+2k-1 & \Rightarrow  & D_j \cdot C_{167} &=&
-j+4l+2k-2 \ge 0, 
\end{array}
\]
and the third case differs for the values of $j$ between $3l-1$ 
and $2l+2k-1$:
\[
\begin{array}{ccccc}
3l-1< j \le l+2k-1 & \Rightarrow  & D_j \cdot C_{167} &=&
3l-1 \ge 0\\
l+2k-1< j \le 2l+2k-1 & \Rightarrow  & D_j \cdot C_{167} &=&
-j+4l+2k-2 \ge 0.
\end{array}
\]
In all three cases, we find $-2j+6l+4k-3 < 0 $ iff $2j \ge 6l+4k-2$ iff 
$j \ge 3l+2k-1 = 2\phi(2)+\phi(3)-1$.
\end{proof}

\begin{thm}
Let $\phi$ be a linear degree function, and 
let  $D_j = jE_0-\sum\limits_{i=1}^7 a_iE_i$, where 
\[
\begin{array}{ccc}
a_1=a_2=a_3 &=& (j-\phi(1)+1)_+\\
a_4=a_5=a_6 &=& (j-2\phi(2)+1)_+\\
a_7 &=& (j-3\phi(3)+1)_+
\end{array}
\] 
Put
\[
\begin{array}{ccc}
t_1 &=& (j-\phi(1)-\phi(2)+2)_+\\
t_2 &=& (j-2\phi(2)-\phi(3)+2)_+
\end{array}
\] 
Then 
$$G= (j-3t_1-3t_2)E_0 -\sum_{i=1}^3(a_i-2t_1-t_2)E_i-\sum_{i=4}^6(a_i-t_1-t_2)E_i-(a_7-3t_2)E_7$$
is a divisor with $h^0(D_j) = h^0(G)$.
If $j \ge \phi(1)+\phi(2)+\phi(3)-2$, then $h^0(G)=0$, and if
$j < \phi(1)+\phi(2)+\phi(3)-2$, then $G$ is an effective, nef 
divisor with 
$$h^0(G)={j+2 \choose 2} -3{j+2-\phi(1) \choose 2}-
3{j+2-2\phi(2) \choose 2}-{j+2-3\phi(3) \choose 2}$$
$$+6 {j+2-\phi(1)-\phi(2) \choose 2}$$
$$+6 {j+2-2\phi(2)-\phi(3) \choose 2}$$
The first row is precisely the dimension when $h^1(D_j)=0$, while the
second row is zero iff $t_1 <2$ and the third row is zero iff $t_2<2$.
\end{thm}
\begin{proof}
Using Lemma 4.2 and Lemma 4.3, $D_j$ can be reduced to a
divisor $G$ which meets the six $-2$ curves non-negatively;
since the $-2$ curves are mutually orthogonal the order of
the reduction is irrelevant. Now, if 
$j \ge \phi(1)+\phi(2)+\phi(3)-2= 3l+3k-2$, then 
$(j-2l-3k+2)_+ = j-2l-3k+2$ and $(j-3l-2k+2)_+ = j-3l-2k+2$. So 
$G \cdot E_0 = -5j+15l+15k-12.$ Our assumption that 
$j \ge 3l+3k-2$ the implies that $G \cdot E_0 < 0$.
In particular, $h^0(G)=h^0(D_j)=0$.

Now assume $j < \phi(1)+\phi(2)+\phi(3)-2= 3l+3k-2$. To show 
that $G$ is effective and nef, we need to show that 
$G\cdot C_{ij} \ge 0$ and $G\cdot E_i \ge 0$. Since the
curves $C_{45},C_{46},C_{56}$ intersect $G$ with the same
multiplicity, to show that $G\cdot C_{ij} \ge 0$ it suffices to
show $G\cdot C_{45} \ge 0$. We have that 
\[
\begin{array}{ccc}
G \cdot C_{45} &= & j -3t_1 - 3t_2 -2(a_4-t_1-t_2)\\
               &= & j -t_1 - t_2 -2a_4\\
               &= & j-(j-2l-3k+2)_+ -(j-3l-2k+2)_+ -2(j-2l-2k+1)_+
\end{array}
\]
There are two cases to consider, depending on the relative values 
of $l,k$. First, suppose $3l+2k-2 \le 2l+3k-2$. 
\[
\begin{array}{ccccc}
0 \le j \le 2l+2k-1 & \Rightarrow  & G \cdot C_{45} &=& j \ge 0\\
2l+2k-1 < j \le 3l+2k-2 & \Rightarrow  & G \cdot C_{45} &=&
-j+4l+4k-2 \ge 0\\
3l+2k-2 < j \le 2l+3k-2 & \Rightarrow  & G \cdot C_{45} &=&
-2j+7l+6k-4 \ge 0\\
2l+3k-2 < j  < 3l+3k-2 & \Rightarrow  & G \cdot C_{45} &=&
-3j+9l+9k-6 \ge 0
\end{array}
\]
The analysis if $2l+3k-2 \le 3l+2k-2$ differs only for $j$ such that:
\[
\begin{array}{ccccc}
2l+2k-1 < j \le 2l+3k-2 & \Rightarrow  & G \cdot C_{45} &=&
-j+4l+4k-2 \ge 0\\
2l+3k-2 < j \le 3l+2k-2 & \Rightarrow  & G \cdot C_{45} &=&
-2j+6l+7k-4 \ge 0
\end{array}
\]
Next, we consider $G\cdot E_i$. Because $E_1, E_2, E_3$
and $E_4, E_5, E_6$ appear with the same multiplicity in $G$, 
it suffices to test $G$ against the classes $E_1, E_4$ and $E_7$.
\[
\begin{array}{ccc}
G \cdot E_1 & = & a_1-2t_1-t_2\\
   & = & (j-\phi(1)+1)_+ -2(j-\phi(1)-\phi(2)+2)_+ -(j-2\phi(2)-\phi(3)+2)_+\\
   & = & (j-l-2k+1)_+ -2(j-2l-3k+2)_+ -(j-3l-2k+2)_+
\end{array}
\]
First, suppose $2l+3k-2 \le 3l+2k-2$. 
\[
\begin{array}{ccccc}
0 \le j \le l+2k-1 & \Rightarrow  & G \cdot E_1 &=& j \ge 0\\
l+2k-1 < j \le 2l+3k-2 & \Rightarrow  & G \cdot  E_1 &=&
j-l-2k+1 \ge 0\\
2l+3k-2< j \le 3l+2k-2 & \Rightarrow  & G \cdot  E_1 &=&
-j+3l+4k-3 \ge 0\\
3l+2k-2 < j < 3l+3k-2 & \Rightarrow  & G \cdot  E_1 &=&
-2j+6l+6k-5 \ge 0
\end{array}
\]
If $3l+2k-2 \le 2l+3k-2$, then the values differ for $j$ such that:
\[
\begin{array}{ccccc}
l+2k-1 < j \le 3l+2k-2 & \Rightarrow  & G \cdot E_1 &=&
j-l-2k+1 \ge 0\\
3l+2k-2 < j \le 2l+3k-2 & \Rightarrow  & G \cdot E_1 &=&
2l-1 \ge 0
\end{array}
\]
In both cases, we see that $G \cdot E_1 < 0$ iff $2j > 6l+6k-5$ iff
$j \ge 3l+3k-2 = \phi(1)+\phi(2)+\phi(3)-2$. The analysis for the
classes $E_4$ and $E_7$ runs along the same lines, with the
same result: $G \cdot E_4 < 0$ iff 
$j \ge 3l+3k-2 = \phi(1)+\phi(2)+\phi(3)-2$ iff $G \cdot E_7 < 0$.

Thus, if $j < \phi(1)+\phi(2)+\phi(3)-2$, then 
Harbourne's results imply that $G$ is an effective, 
nef divisor, with $h^1(G)=h^2(G)=0$, and $h^0(G)=h^0(D_j)$. 
This reduces the computation of $h^0(D_j)$ to a
simple numerical calculation, the point being that for $G$, 
$h^1(G)=h^2(G)=0$, so:
$$h^0(G)= \chi(G)= \frac{G^2-G K_X}{2}+1.$$ 
When $t_1=1$, the second row of the formula does not contribute 
to $h^0(G)$ because if $t_1=1$, then 
$D'_j = D_j - C_{124}-C_{135}-C_{236}$ and a calculation shows that 
$$D^2_j-D_j \cdot K =D_j'^2-D_j' \cdot K.$$
This stems from the choice of a degree function which is linear; a 
similar observation pertains when $t_2=1$.
\end{proof}
In \cite{h2}, Harbourne works over an algebraically closed field $\mathbb{K}$, 
but as long as $\mathbb{K}$ contains the 
coordinates of the points being blown up, this hypothesis is
unnecessary. This is because the rank of the matrix 
which computes the dimension of the space of global sections 
is the same over $\mathbb{K}$ or $\overline{\mathbb{K}}$. 
However, to apply the results of Section 2, we need to assume 
that $char(\mathbb{K})=0$ or $char(\mathbb{K}) > \sum_{i=1}^3\phi(i)-2$. 
In this case, the theory of inverse systems tells us that 
$\mbox{dim}_\mathbb{K}H^0(D_j)=\mbox{dim}_\mathbb{K}(R/J_\phi)_j$.
A simple computation translates this data 
into the Hilbert series:
\begin{corollary} For $n=3$ and $\phi$ linear, if $char(\mathbb{K})
  > \phi(1)+\phi(2)+\phi(3)-2$ or $char(\mathbb{K})=0$, then $P(R/J_\phi, t)=$
$$\frac{1-3t^{\phi(1)}-3t^{2 \phi(2)}-t^{3\phi(3)}+6t^{\phi(1)+\phi(2)}+6t^{2\phi(2)+\phi(3)}-6t^{\phi(1)+\phi(2)+\phi(3)}}{(1-t)^3}.$$
\end{corollary}
Computing the Hilbert series of $R/I_\phi$ via the Taylor resolution 
mentioned earlier and comparing with the result above proves 
Theorem 1.1 for a linear degree function. 
While it may be possible to obtain a better 
lower bound on $char(\mathbb{K})$, Example 3.2 shows that this 
bound is sometimes tight.

\section{Almost linear degree functions}
A degree function $\phi$ (with values in $\mathbb{N}$) 
is called {\it almost linear} if there is a $k \in \mathbb{N}$ such that 
$$\phi(r)-\phi(r+1)=k\mbox{ or }k+1\mbox{ }\forall r \in \{1,\ldots,n-1\}.$$
When $n=3$, an almost linear degree function satisfies the same key
property satisfied by a linear degree function: if $D_j$ is effective, 
then it can be reduced to an effective, nef divisor using only the 
classes $C_{ijk}$. 
For an almost linear degree function, the Hilbert series of $R/J_\phi$
will have two possible forms, depending on the difference vector $\Delta = (\phi(1)-\phi(2),\phi(2)-\phi(3))$. 
\begin{lemma}
Let $\phi$ be an almost linear degree function.
Then $D_j$ meets $C_{124},C_{135},C_{236}$
negatively iff $j \ge \phi(1)+\phi(2)-1$. Put 
$$t_1= (j-\phi(1)-\phi(2)+2)_+$$
Then $$D_j'=(j-3t_1)E_0-\sum_{i=1}^3(a_i-2t_1)E_i-\sum_{i=4}^6(a_i-t_1)E_i-a_7E_7$$
meets $C_{124},C_{135},C_{236}$ non-negatively, and has $h^0(D_j)=h^0(D_j')$.
\end{lemma}
\begin{proof}
We have to consider two cases. First, suppose $\Delta = (k+1,k)$. 
For some $l \in \mathbb{N}$ we have 
$\phi(3) = l$, so $\phi(2) = l+k$ and $\phi(1) = l+2k+1$.
\[
\begin{array}{ccc}
D_j \cdot C_{124} &=& j-2(j-\phi(1)+1)_+ -(j-2\phi(2)+1)_+\\
                  & = & j-2(j-l-2k)_+ -(j-2l-2k+1)_+
\end{array}
\]

\[
\begin{array}{ccccc}
0 \le j \le l+2k & \Rightarrow  & D_j \cdot C_{124} &=& j \ge 0\\
l+2k < j \le 2l+2k-1 & \Rightarrow  & D_j \cdot C_{124} &=&
-j+2l+4k \ge 0\\
2l+2k-1 < j & \Rightarrow  & D_j \cdot C_{124} &=&
-2j+4l+6k-1
\end{array}
\]
So $D_j \cdot C_{124} <0$ iff $2j \ge 4l+6k$ iff $j \ge 2l+3k =
\phi(1)+\phi(2)-1$. Now suppose $\Delta = (k,k+1)$, and let
$\phi(3) = l$, so $\phi(2) = l+k+1$ and $\phi(1) = l+2k+1$.
\[
\begin{array}{ccc}
D_j \cdot C_{124} &=& j-2(j-\phi(1)+1)_+ -(j-2\phi(2)+1)_+\\
                  & = & j-2(j-l-2k)_+ -(j-2l-2k-1)_+
\end{array}
\]
Proceeding as in the previous case, we find 
$D_j \cdot C_{124} <0$ iff $2j \ge 4l+6k+2$ iff $j \ge 2l+3k+1 =
\phi(1)+\phi(2)-1$.
\end{proof}

\begin{thm}
Let $\phi$ be an almost linear degree function; 
define $a_i$ and $D_j$ as in the 
previous section. Put
$$t_1 = (j-\phi(1)-\phi(2)+2)_+,\mbox{ }t_2 = (j-2\phi(2)-\phi(3)+2)_+$$
Then 
$$G= (j-3t_1-3t_2)E_0 -\sum_{i=1}^3(a_i-2t_1-t_2)E_i-\sum_{i=4}^6(a_i-t_1-t_2)E_i-(a_7-3t_2)E_7$$
is a divisor with $h^0(D_j) = h^0(G)$.
If $j \ge \phi(1)+\phi(2)+\phi(3)-2$, then $h^0(G)=0$, and if
$j < \phi(1)+\phi(2)+\phi(3)-2$, then $G$ is effective and nef, with $h^0(D_j) = h^0(G) =$
$${j+2 \choose 2} -3{j+2-\phi(1) \choose 2}-
3{j+2-2\phi(2) \choose 2}-{j+2-3\phi(3) \choose 2}$$
$$+6{j+2-\phi(1)-\phi(2) \choose 2}$$
$$+3{j+2-\phi(3)-2\phi(2) \choose 2}+3{j+2-\phi(1)-2\phi(3) \choose 2}$$
The first row is precisely the dimension when $h^1(D_j)=0$, while the
second row is zero iff $t_1 <2$ and the third row is zero iff $t_2<1$.
\end{thm}
\begin{proof}
First, suppose $\Delta = (k,k+1)$, so for some $l \in \mathbb{N}$ we can write 
$\phi(3) = l$, $\phi(2) = l+k+1$ and $\phi(1) = l+2k+1$.
\[
\begin{array}{ccc}
D_j \cdot C_{167} &= & j-(j-\phi(1)+1)_+ -(j-2\phi(2)+1)_+ -(j-3\phi(3)+1)_+\\
 &= & j-(j-l-2k)_+ -(j-2l-2k-1)_+ -(j-3l+1)_+
\end{array}
\]
We have to analyze three possible cases
\[
\begin{array}{ccccc}
l+2k &\le& 2l+2k+1& \le& 3l-1\\
l+2k &\le& 3l -1 &\le & 2l+2k+1\\
3l-1 &\le & l+2k &\le& 2l+2k+1
\end{array}
\]
Working through the possibilities, we find $D_j \cdot C_{167} < 0$ iff
$j \ge 3l+2k+1 = \phi(1)+2\phi(3)$. So removing 
$t_2 = (j - \phi(1)-2\phi(3)+1)_+$ copies of $C_{167}+C_{257}+C_{347}$
yields a divisor which meets each of the classes 
$\{ C_{167},C_{257},C_{347} \}$ non-negatively; and applying Lemma 5.1 
results in a divisor $G$ which meets all the $-2$-curves
non-negatively. If $j \ge \phi(1)+\phi(2)+\phi(3)-2= 3l+3k$, then 
$G \cdot E_0 = -5j+15l+15k$. So if $j > \phi(1)+\phi(2)+\phi(3)-2$, 
$G \cdot E_0 <0$ and $h^0(G)=0$. If $j = \phi(1)+\phi(2)+\phi(3)-2$,
then $G = -E_4-E_5-E_6+E_7$. Since $G \cdot C_{45} <0$, the
algorithm tells us to further reduce $G$ by removing $C_{45}$, which
leaves the divisor $-E_0-E_6+E_7$, hence $h^0(G)=0$.
 
Now suppose $j < \phi(1)+\phi(2)+\phi(3)-2= 3l+3k$. As in the
proof of Theorem 4.4,  to show that $G$ is effective and nef it 
suffices to show that $G$ meets the
classes $\{C_{45},E_1,E_4, E_7 \}$ non-negatively. 
\[
\begin{array}{ccc}
G \cdot C_{45} &= & j -3t_1 - 3t_2 -2(a_4-t_1-t_2)\\
               &= & j -t_1 - t_2 -2a_4\\
               &= & j-(j-2l-3k)_+ -(j-3l-2k)_+ -2(j-2l-2k-1)_+
\end{array}
\]
There are, again, two cases, depending on if 
$3l+2k \le 2l+3k$ or $2l+3k \le 3l+2k$. First, suppose $3l+2k \le 2l+3k$. 
\[
\begin{array}{ccccc}
0 \le j \le 2l+2k+1 & \Rightarrow  & G \cdot C_{45} &=& j \ge 0\\
2l+2k+1 < j \le 3l+2k & \Rightarrow  & G \cdot C_{45} &=&
-j+4l+4k+2 \ge 0\\
3l+2k < j \le 2l+3k & \Rightarrow  & G \cdot C_{45} &=&
-2j+7l+6k+2 \ge 0\\
2l+3k < j  < 3l+3k & \Rightarrow  & G \cdot C_{45} &=&
-3j+9l+9k+2 \ge 0
\end{array}
\]
If $2l+3k \le 3l+2k$, then the values differ for $j$ such that:
\[
\begin{array}{ccccc}
2l+2k+1 < j \le 2l+3k & \Rightarrow  & G \cdot C_{45} &=&
-j+4l+4k+2 \ge 0\\
2l+3k < j \le 3l+2k & \Rightarrow  & G \cdot C_{45} &=&
-2j+6l+7k+2 \ge 0 
\end{array}
\]
So $G \cdot C_{45} \ge 0$ when $j < 3l+3k$. Next, we test 
$G$ against $E_1, E_4$ and $E_7$.
\[
\begin{array}{ccc}
G \cdot E_1 & = & a_1-2t_1-t_2\\
   & = & (j-\phi(1)+1)_+ -2(j-\phi(1)-\phi(2)+2)_+ -(j-\phi(1)-2\phi(3)+1)_+\\
   & = & (j-l-2k)_+ -2(j-2l-3k)_+ -(j-3l-2k)_+
\end{array}
\]
First, suppose $2l+3k \le 3l+2k$. 
\[
\begin{array}{ccccc}
0 \le j \le l+2k & \Rightarrow  & G \cdot E_1 &=& j \ge 0\\
l+2k < j \le 2l+3k & \Rightarrow  & G \cdot  E_1 &=&
j-l-2k \ge 0\\
2l+3k< j \le 3l+2k & \Rightarrow  & G \cdot  E_1 &=&
-j+3l+4k \ge 0\\
3l+2k < j < 3l+3k & \Rightarrow  & G \cdot  E_1 &=&
-2j+6l+6k \ge 0
\end{array}
\]
If $3l+2k \le 2l+3k$, then the values differ for $j$ such that:
\[
\begin{array}{ccccc}
l+2k < j \le 3l+2k & \Rightarrow  & G \cdot E_1 &=&
j-l-2k \ge 0\\
3l+2k < j \le 2l+3k & \Rightarrow  & G \cdot E_1 &=&
2l \ge 0
\end{array}
\]
In both cases, we see that $G \cdot E_1 < 0$ if $2j > 6l+6k$ if
$j \ge 3l+3k +1 = \phi(1)+\phi(2)+\phi(3)-1$. A similar analysis 
shows that $G \cdot E_4 < 0$ iff $j \ge \phi(1)+\phi(2)+\phi(3)-2$,
and $G \cdot E_7 < 0$ iff $j \ge \phi(1)+\phi(2)+\phi(3)-1$.
This concludes the proof if $\Delta = (k,k+1)$. 

A similar argument yields the Hilbert function if $\Delta = (k+1,k)$;
when expressed in terms of $k,l$ it differs from the formula for
$\Delta = (k,k+1)$. But because
\[
\begin{array}{ccc}
\Delta = (k,k+1)& \Rightarrow  & 2\phi(2)+\phi(3) =  3l +2k+ 2\\
\Delta = (k+1,k)& \Rightarrow  & 2\phi(2)+\phi(3) =  3l +2k,
\end{array}
\]
writing the Hilbert function in terms of $\phi$ yields a formula
independent of $\Delta$.
\end{proof}
\begin{corollary} For $n=3$ and $\phi$ almost linear, if $char(\mathbb{K})
  > \phi(1)+\phi(2)+\phi(3)-2$ or $char(\mathbb{K})=0$, then $P(R/J_\phi, t)=$
\begin{small}
\[
\frac{1-3t^{\phi(1)}-3t^{2 \phi(2)}-t^{3\phi(3)}+
3t^{2\phi(2)+\phi(3)}+3t^{\phi(1)+2\phi(3)}+6t^{\phi(1)+\phi(2)}
-6t^{\phi(1)+\phi(2)+\phi(3)}}{(1-t)^3}.
\]
\end{small}
\end{corollary}
\begin{example}
Let $char(\mathbb{K}) = 0$. The almost linear degree function 
\[
\phi(1)=16, \phi(2)=12,\phi(3)=7,
\]
has $\Delta = (4,5)$ and associated Hilbert series 
\[
\frac{1-3t^{16}-t^{21}-3t^{24}+6t^{28}+3t^{30}+3t^{31}-6t^{35}}{(1-t)^3}.
\]
The almost linear degree function 
\[
\phi(1)=20, \phi(2)=16, \phi(3)=13,
\]
has $\Delta = (4,3)$ and associated Hilbert series
\[
\frac{1-3t^{20}-3t^{32}+6t^{36}-t^{39}+3t^{45}+3t^{46}-6t^{49}}{(1-t)^3}.
\]
\end{example}
Postnikov and Shapiro note that when the
degree function is not linear or almost linear, they were unable to find
a case where the equality of Theorem 1.1 held. In general, Harbourne's
methods can be applied to determine the Hilbert series of any ideal
generated by powers of linear forms in three variables, as long as
there are at most eight generators (and in certain cases, more). 

\section{Examples, conjectures and a counterexample}
We assume throughout this section that $char(\mathbb{K})
  > \phi(1)+\phi(2)+\phi(3)-2$ or $char(\mathbb{K})=0$.
We use $G_\phi$ to denote an ideal of generic 
forms generated in the same degrees as $J_\phi$.
\begin{example}
We revisit Example 3.3. For this example we have:
$$P(R/G_\phi,t) = P(R/J_\phi,t) = P(R/I_\phi,t),$$
the Hilbert series for all three is
$$1 + 3t + 6t^2  + 9t^3  + 12t^4  + 12t^5  + 6t^6.$$
\end{example}
\begin{example}
The degree function $\phi(1)=8, \phi(2)=5, \phi(3)=3$ is almost
 linear, but not linear. 
We illustrate the computation of $h^0(D_{13})$. Since $t_1=2$ we first
reduce to $D_{13}-2(C_{124}+C_{135}+C_{236})$, yielding 
$$D'= 7E_0- \sum_{i=1}^6 2E_i-5E_7.$$
Finally, we check that $t_2=1$ (this is the almost linear $t_2$), so 
$$G=4E_0- \sum_{i=1}^6 E_i-2E_7.$$
We find $$h^0(G)= \frac{G^2-G K_X}{2}+1=6.$$
Of course, Corollary 5.3 tells us that 
$$P(R/J_\phi,t) = \frac{1-t^6-3t^{8}-3t^{10}+9t^{13}+3t^{14}-6t^{16}}{(1-t)^3},$$
so this is the expected value.
\end{example}
We close the paper with a conjecture from \cite{ssp} about 
the minimal free resolution of $J_\phi$ when $\phi$ is linear. 
Roughly speaking, the conjecture is that 
the resolutions are as simple as can be expected. In \cite{ps},
Postnikov and Shapiro generalize this conjecture to a 
broader class of families.

\begin{conj}
For $n=3$ and $\phi$ linear, the minimal free 
resolution of $J_\phi$ is given by:
$$0 \longrightarrow R^6(-\sum\limits_{i=1}^3\phi(i)) \longrightarrow
\begin{array}{c}
 R^6(-2\phi(2)-\phi(3))\\
\oplus\\
R^6(-\phi(1)-\phi(2)) \\
\end{array}
\longrightarrow
\begin{array}{c}
R^3(-\phi(1))\\ 
\oplus\\
R^3(-2\phi(2)) \\
 \oplus \\
R(-3\phi(3))\\
\end{array}
\longrightarrow
J_\phi \longrightarrow 0.$$
\end{conj}

A similar pair of conjectures can be made for the almost linear case;
when $\Delta = (k+1,k)$ the conjecture is the natural analog of
Conjecture 6.3. However, interesting things happen when $\Delta =
(k,k+1)$.
\begin{example}
For the almost linear degree function
$$\phi(1)=8, \phi(2)=5, \phi(3)=1,$$
the minimal free resolution of $J_\phi$ is given by:
$$0 \rightarrow R^6(-14) \rightarrow
 R^6(-13) \oplus R^3(-11)  \rightarrow R^3(-8) \oplus R(-3)
 \rightarrow J_\phi \rightarrow 0$$
and the minimal free resolution of $I_\phi$ is given by:
$$0 \rightarrow R^6(-14) \rightarrow
\begin{array}{c}
 R^6(-13)\\ 
\oplus\\
R^3(-11) \\
 \oplus \\
 R^3(-10)\\
\end{array}
\rightarrow
\begin{array}{c}
 R^3(-10)\\ 
\oplus\\
R^3(-8) \\
 \oplus \\
R(-3)\\
\end{array}
 \rightarrow I_\phi \rightarrow 0.$$
So while $P(R/J_\phi,t)=P(R/I_\phi,t)$ as required by Corollary
 5.3, the minimal free resolutions differ. In particular, this gives
a counterexample to Conjecture 6.10 of \cite{ps} (see also remarks 
after Corollary 12.2).
\end{example}

\end{document}